# $M|D|\infty$ QUEUE BUSY PERIOD DISTRIBUTION

Prof. Dr. **MANUEL ALBERTO M. FERREIRA**

Instituto Universitário de Lisboa (ISCTE – IUL), BRU - IUL, Lisboa, Portugal

manuel.ferreira@iscte.pt

**ABSTRACT**

The M|D|∞ queue system busy period distribution is intensively studied because, due to its probability density function quite easy interpretation, it may serve as a clue to interpret M|G|∞ systems, with other service time's distributions, busy period distributions. In this work the M|D|∞ queue system busy period distribution parameters are exactly computed and its probability density function is interpreted and explained. The only problem is how to compute the distribution function, for which an algorithm is presented and exemplified its application.

**Keywords**: M|D|∞, busy period, probability density function, distribution function.

## 1. INTRODUCTION

A queue system busy period is a period that begins when a customer arrives at the system finding it empty, ends when a customer abandons the system letting it empty and, in it, there is always at least one customer present.

In the M|G|∞ queue system the customers arrive according to a Poisson process at rate $\lambda$, receive a service which time length is a positive random variable with distribution function $G(.)$ and mean $\alpha$ and, when they arrive, each one finds immediately an available server. Each customer service is independent from the other customers' services and from the arrivals process. The traffic intensity is $\rho = \lambda\alpha$.

It is intended, in this work, to study the M|D|∞ queue system busy period distribution. In this system the service is constant with value α, in other words: **D**eterministic with value α.

Call $B$ the busy period length random variable and $b(t), B(t)$ and $\bar{B}(s)$ its probability density function, distribution function and Laplace transform respectively.

$\bar{B}(s)$ may be described in terms of $G(.)$ and $\lambda$ as

$$\bar{B}(s) = 1 + \lambda^{-1}\left(s - \frac{1}{\int_0^\infty e^{-st-\lambda\int_0^t[1-G(v)]dv}dt}\right) \qquad (1.1),$$

see (1).

This transform inversion is not easy, apart from some cases (see 2, 3 and 4).

In the M|D|∞ system, the probability density function obtained inverting (1.1), although complex is easy to interpret since, in this system, the order for which the customers abandon is the same for which they arrive at the system. And it is possible to obtain exact expressions for the $B$ moments The same does not happen for $B(t)$.

In section 2 is studied the M|D|∞ system $b(t)$. In section 3 the M|D|∞ system $B$ moments. In section 4 the analytical results possible for the M|D|∞ system $B(t)$ are presented. Finally, in section 5, an algorithm that allows computing approximately $B(t)$ for the M|D|∞ system is presented.

## 2. *B* PROBABILITY DENSITY FUNCTION FOR THE M|D|∞ SYSTEM

For the M|D|∞ queue system (1.1) has the form

$$\bar{B}(s) = 1 + \lambda^{-1}\left(s - \frac{(s+\lambda)s}{\lambda e^{-(s+\lambda)\alpha} + s}\right) \quad (2.1)$$

obtaining, by Laplace transform inversion, see (5),

$$b(t) = \sum_{n=0}^{\infty} \left(\frac{d}{dt}\frac{c(t)}{e^{-\rho}}\right) * \left(\frac{d}{dt}\frac{1-d(t)}{1-e^{-\rho}}\right)^{*n} e^{-\rho}(1 - e^{-\rho})^n \quad (2.2)$$

where

- * is the convolution operator,

- $\frac{c(t)}{e^{-\rho}} = \begin{cases} 0, t < \alpha \\ 1, t \geq \alpha \end{cases} = G(t),$

- $\frac{1-d(t)}{1-e^{-\rho}} = \begin{cases} \frac{1-e^{-\lambda t}}{1-e^{-\rho}}, t < \alpha \\ 1, t \geq \alpha \end{cases}.$

Note that the M|D|∞ queue system output in equilibrium is a Poisson process with parameter $\lambda$, see (6). So (2.2) may be written as

$$b(t) = \sum_{n=0}^{\infty} g(t) * \left(\frac{d}{dt}P(t_o \leq t | t_o \leq \alpha)\right)^{*n} e^{-\rho}(1 - e^{-\rho})^n \quad (2.3)$$

where $g(\cdot)$ is the service probability density function and $t_o$ is the time that occurs between two consecutive outputs from the queue.

Then

- $g(t)$ corresponds to the first customer contribution to the busy period length

- $\frac{d}{dt}P(t_o \leq t | t_o \leq \alpha)$ corresponds to the $n^{th}$ customer contribution to the busy period length, $n \geq 2$, $e^{-\rho}(1 - e^{-\rho})^n$

$-e^{-\rho}(1-e^{-\rho})^n$ is the probability that the busy period ends in the $n+1$ output, $n = 0,1,2,...$

and, so,

$$-g(t) * \left(\frac{d}{dt}P(t_o \leq t | t_o \leq \alpha)\right)^{*n} \text{ is the probability density function of a}$$

busy period with $n+1$ customers, $n = 0,1,2,...$ . This period will have a length lesser or equal than $(n+1)\alpha, n = 0,1,2,...$ .

### 3. $B$ MOMENTS FOR THE M|D|∞ SYSTEM

The expression (1.1) is equivalent to $(\bar{B}(s) - 1)(C(s) - 1) = \lambda^{-1}sC(s)$ where $C(s) = \int_0^\infty e^{-st - \lambda \int_0^t [1-G(v)]dv} \lambda(1-G(t))dt$. Computing the $n^{\text{th}}$ derivative, using Leibnitz's formula and making $s = 0$, it is obtained:

$$E[B^n] = (-1)^{n+1}\left\{\frac{e^\rho}{\lambda}n\, C^{(n-1)}(0) - e^\rho \sum_{p=1}^{n-1}(-1)^{n-p}\binom{n}{p}E[B^{n-p}]C^{(p)}(0)\right\},$$

$$n = 1,2,... \quad (3.1)$$

where

$$C^{(n)}(0) = \int_0^\infty (-t)^n e^{-\lambda \int_0^t [1-G(v)]dv} \lambda(1-G(t))dt, n = 01,2,... \quad (3.2).$$

For the M|D|∞ system:

$$C^{(0)}(0) = 1 - e^{-\rho}$$
$$C^{(n)}(0) = -e^{-\rho}(-\alpha)^n - \frac{n}{\lambda}C^{(n-1)}(0), \quad n = 1,2,... \quad (3.3)$$

and it is possible to compute exact values for the any order around the origin moments, not being necessary to use numerical methods, as it happens for the generality of the M|G|∞ systems, see (7). Then

$-n = 1$

$$E[B] = \frac{e^\rho - 1}{\lambda} \quad (3.4)$$

as in any other M|G|∞ system,

$-n = 2$

$$E[B^2] = \frac{2e^\rho(e^\rho - 1 - \rho)}{\lambda^2} \quad (3.5).$$

From (3.5) and (3.4) it results

$$VAR[B] = \frac{2e^\rho(e^\rho - 1 - \rho) - (e^\rho - 1)^2}{\lambda^2} \quad (3.6)$$

and from (3.6) and (3.4)

$$\gamma_B^2 = \frac{2e^\rho(e^\rho - 1 - \rho)}{(e^\rho - 1)^2} - 1 \quad (3.7)$$

being $\gamma_B$ the $B$ variation coefficient.

From (3.7) it is possible to conclude that

- $\gamma_B$ depends only on $\rho$,

- $\gamma_B < 1$,

- $\lim_{\rho \to \infty} \gamma_B = 1$ (note that it is the exponential variation coefficient value).

In Table 1 are shown some values of $\gamma_B$, with eight significant figures, computed giving to $\rho$ values between *0.5* and *100*. In practice, after $\rho = 10$,

**Table 1**. $\gamma_B$ values for the M|D|∞ system

| $\rho$ | $\gamma_B$ |
|---|---|
| .5 | .40655883 |
| 1 | .56798436 |
| 10 | .99959129 |
| 20 | .99999999 |
| 50 | .99999999 |
| 100 | .99999999 |

they do not distinguish from *1*.

Calling $\beta_1$ the $B$ asymmetry coefficient, see for instance (8, 9), its value is given by

$$\beta_1 = \frac{(2e^{3\rho} - 6(1-\rho)e^{2\rho} - 6(1+\rho)e^\rho - 12\rho + 10)^2}{(e^{2\rho} - 2\rho e^\rho - 1)^3}, \rho \neq 0 \quad (3.8).$$

Evidently, $\lim_{\rho \to \infty} \beta_1 = 4$ (note that it is the exponential asymmetry coefficient value).

Computing some values for $\beta_1$ in the same circumstances of the ones for $\gamma_B$, see Table 2, it is evident that after $\rho = 10$ $\beta_1$ in practice does not distinguish from 4.

**Table 2.** $\beta_1$ values for the M|D|∞ system

| $\rho$ | $\beta_1$ |
|---|---|
| .5 | 6.0360869 |
| 1 | 4.5899937 |
| 10 | 4.0000001 |
| 20 | 4.0000000 |
| 50 | 4.0000000 |
| 100 | 4.0000000 |

and it is always $\beta_1 > 0$ what indicates a distribution right asymmetric.

Calling $\beta_2$ the B kurtosis, see still (8, 9), its value is given by

$$\beta_2 = (9e^{4\rho} - 24e^{3\rho} + 96e^{2\rho} + 36e^{\rho} - 117 + 72\rho e^{\rho} - 60\rho e^{3\rho} - 96\rho e^{2\rho}$$
$$- 72\rho^2 e^{\rho} + 24\rho^2 e^{2\rho} - 4\rho^3 - 36\rho^2 - 160\rho)/(e^{2\rho} - 2\rho e^{\rho} - 1)^2,$$
$$\rho \neq 0 \qquad (3.9).$$

As it is evident, $\lim_{\rho \to \infty} \beta_2 = 9$ (note that it is the exponential kurtosis value). Computing some values for $\beta_2$ in the same circumstances of the ones for $\gamma_B$, see Table 3, it is evident that after $\rho = 10$ $\beta_2$ in practice does not distinguish from 9.

**Table 3.** $\beta_2$ values for the M|D|∞ system

| $\rho$ | $\beta_2$ |
|---|---|
| .5 | 11.142336 |
| 1 | 9.6137084 |
| 10 | 9.0000000 |
| 20 | 9.0000000 |
| 50 | 9.0000000 |
| 100 | 9.0000000 |

It is always $\beta_2 > 3$ what indicates a distribution leptokurtic, that is: with a very sharp peak.

### 4. THE DISTRIBUTION FUNCTION OF B FOR THE M|D|∞ SYSTEM

After (2.1) it may be obtained

$$B(t) = 1 - \lambda^{-1} \sum_{n=1}^{\infty} \left( \frac{d}{dt} \frac{1-d(t)}{1-e^{-\rho}} \right)^{*n} (1-e^{-\rho})^n \quad (4.1).$$

For details see (5). But, for α great enough,

$$\frac{d}{dt} \frac{1-d(t)}{1-e^{-\rho}} \cong \lambda e^{-\lambda t}, t \geq 0 \quad (4.2).$$

So, under these conditions,

$$B(t) \cong 1 - e^{-\lambda e^{-\rho} t}, t \geq 0 \quad (4.3)$$

that is: B is approximately exponential in coherence with what it was seen in section 3 for $\gamma_B$, $\beta_1$ and $\beta_2$; confer with (10). And, see still (5),

$$e^{-\rho} G(t) \leq B(t) \leq G(t), t \geq 0 \quad (4.4).$$

So, for values of $\rho$ close from 0, $B(t) \cong G(t)$. And, being $t_p$ the quantile corresponding to the probability p, for B after (4.4), for the M|D|∞ system it is concluded that

$$\begin{array}{l} t_p = \alpha, \ p \leq e^{-\rho} \\ t_p \geq \alpha \end{array} \quad (4.5).$$

Finally, using the Chebyshev inequality, see again (5)

$$B(t) \geq 1 - \frac{e^{2\rho} - 2\rho e^{\rho} - 1}{(1+\lambda t + e^{\rho})^2} \text{ since } t > \lambda^{-1} \left[ e^{\rho} - 1 + \text{máx} \left\{ e^{\rho} - 1; \sqrt{e^{2\rho} - 2\rho e^{\rho} - 1} \right\} \right] \quad (4.6).$$

## 5. ALGORITHM FOR THE APPROXIMATE COMPUTATION OF *B(t)* FOR THE M|D|∞ SYSTEM

As there is no practical expression to calculate the M|D|∞ busy period distribution function, except for very low or very high values of $\rho$ as it was seen in section 4, it will be presented in this section a method to calculate $B(t)$, based in an algorithm of Platzman, Ammons and Bartholdi III, see (11).

Be X a random variable and A a given number. Platzman, Ammons and Bartholdi III present an algorithm to compute an approximate value of $P[X > A]$ that demands the knowledge of the X Laplace transform, $L(s)$, in simple form. So, in principle, this algorithm may be applied to calculate $B(t)$ for any M|G|∞ system since $\bar{B}(s)$ is known, expression (1.1). But, indeed, this does not happen because only for the cases referred in section 1, see again (2, 3 and 4), for which the $\bar{B}(s)$ calculation does not present any problem, and constant service times: expression (2.1), $\bar{B}(s)$ has a simple form.

It is generally said that:

-An algorithm is "accurate" if it looks for solving a problem "close" to the one that is supposed to solve,

-An algorithm is "precise" if it gets a solution "close" to the one of the problem that it is trying to solve.

More concretely, being $\Delta A$ ($\Delta A > 0$) the accuracy and $\Delta p$ $\left(0 < \Delta p < \frac{1}{2}\right)$ the precision required, the approximation $\tau$ of $P[X > A]$ must satisfy the condition

$$P[X \geq A + \Delta A] - \Delta p \leq \tau \leq P[X > A - \Delta A] + \Delta p \quad (5.1).$$

Platzman, Ammons and Bartholdi III suggest doing

$$\tau = \frac{U-A+\Delta A}{U-L+2\Delta A} + \sum_{n=1}^{N} \frac{\alpha^{n^2}}{\pi n} im\{(\beta^n - \gamma^n)L(j\omega n)\} \quad (5.2)$$

where

-$K = \log \frac{2}{\Delta p}$,

-$D = \frac{\Delta A}{\sqrt{2K}}$,

-$\omega = \frac{2\pi}{U-L+2\Delta A}$,

-$N = \left[\frac{2K}{\omega \Delta A}\right]$, being $[\cdot]$ the characteristic of a real number,

-$\alpha = e^{-D^2 \frac{\omega^2}{2}}$,

-$\beta = e^{j(U+\Delta A)\omega}$

-$\gamma = e^{jA\omega}$,

-$U$ and $L$ are numbers such that $1 - P[L \leq X \leq U] \ll \Delta p$,

-$j = \sqrt{-1}$ and

-$im(\cdot)$ designates the imaginary part of a complex number and demonstrate that the approximation so defined fulfills the condition (5.1).

Generally, this method may be described as follows:

-To guarantee a quick execution only $N$ values of the transform are computed. Those values are carefully selected in order to supply as much information as possible. Then the exact value of the tail corresponding to the most smooth distribution function, which Laplace transform passes by those N points, is computed.

It is expected that a method of this kind behaves at least as well as any other method that computes N transform values, and quicklier than one that has to compute it more times.

Platzman, Ammons and Bartholdi III show also that to compute a tail after a transform is a $\#P-$ hard problem. This is indicative of the computational work demanded because to resolve a $\#P-$ hard problem, even with a certain guarantee of approximation, requires an increase of computation that grows exponentially with the problem dimension. Note also that the algorithm does not give a solution of the original problem, but from an approximation defined by $\Delta A$ and $\Delta p$.

Note that in the definition of error used, $\Delta A$ refers to a disturbance of the parameter $A$ while the error most common definition refers only to a disturbance $\Delta p$ of the result.

It is possible to use this algorithm in the computation of the M|D|∞ queue busy period distribution function because, in this case, the Laplace transform has a simple expression, confer (2.1).

To do so put:

-$L = a$ (instead of $\alpha$ for obvious reasons),

-$U = \lambda^{-1}\left(e^\rho - 1 + \sqrt{\frac{e^{2\rho}-2\rho e^\rho-1}{\Delta p}10^l}\right), l = 1,2,\ldots$ since $1 - P[L \leq X \leq U] = 1 - P[a \leq X \leq U] = 1 - P[0 \leq X \leq U]$ having to be $B(U) > 1 - 10^l \Delta p$ and this happens since $\frac{e^{2\rho}-2\rho e^\rho-1}{(1+\lambda U-e^\rho)^2} = 10^{-l}\Delta p \Leftrightarrow (1+\lambda U - e^\rho)^2 = \frac{e^{2\rho}-2\rho e^\rho-1}{\Delta p}10^l$ what leads to the value presented, remember (4.6),

-$A = t$(time),

-the obtained values $B_c(t)$ are given by $1 - \tau$.

Taking $l = 3$ a FORTRAN program was built to implement the algorithm, see (12). The values of $a, A = t, \Delta A = \Delta t$ and $\Delta p$ must be indicated.

In Table 4 are presented results of the application of this algorithm to the computation of the M|D|∞ queue busy period distribution function in the following cases:

- $\alpha = .1$ and $\lambda = 1$,
- $\alpha = 1$ and $\lambda = 1$,
- $\alpha = 3$ and $\lambda = 1$.

**Table 4**. Busy period distribution function computation for the M|D|∞ system

- $\alpha =.1; \lambda = 1; \rho =.1$

| | | $\Delta t = .001; \Delta p = .001$ | |
|---|---|---|---|
| t | $B_1(t)$ | $B_2(t)$ | $B_c(t)$ |
| .11 | -14.805062 | .904837 | .94131 |
| .15 | .816597 | .904837 | .950782 |
| .2 | 959013 | .904837 | .996209 |
| .25 | .982428 | .904837 | .999575 |

|  | Exact | Computed after $B_c(t)$ with $B_c(.1) = .904837$ | Error |
|---|---|---|---|
| $E[B]$ | .105170918 | .1049714128 | 0,2% |
| $VAR[B]$ | .0003685744 | .00031661238 | 14% |

- $\alpha = 1; \lambda = 1; \rho = 1$

$\Delta t = .1; \Delta p = .001$

| $t$ | $B_1(t)$ | $B_2(t)$ | $B_c(t)$ |
|---|---|---|---|
| .2 | -11.001397 | .367879 | .741497 |
| .3 | .420202 | .367879 | .907228 |
| .4 | .817048 | .367879 | .969885 |
| .5 | .911558 | .367879 | .992784 |

|  | Exact | Computed after $B_c(t)$ with $B_c(1) = .904837$ | Error |
|---|---|---|---|
| $E[B]$ | 1.718281828 | 1.6649785 | 3% |
| $VAR[B]$ | .9524924414 | .70343785 | 26% |

- $\alpha = 3; \lambda = 1; \rho = 3$

$\Delta t = .5; \Delta p = .01$

| $t$ | $B_1(t)$ | $B_2(t)$ | $B_c(t)$ |
|---|---|---|---|
| 4 | -.238790 | .0497871 | .099527 |
| 5 | -.420929 | .0497871 | .148885 |
| 6 | -.646402 | .0497871 | .198405 |
| 7 | -.930133 | .0497871 | .244893 |
| 8 | -1.294064 | .0497871 | .288204 |
| 9 | -1.771539 | .0497871 | .329391 |
| 10 | -2.415214 | .0497871 | .368208 |
| 15 | -15.889655 | .0497871 | .530699 |
| 20 | -336.121704 | .0497871 | .65134 |
| 25 | -7.0691347 | .0497871 | .740937 |
| 30 | -1.366543 | .0497871 | .807469 |
| 35 | -.133102 | .0497871 | .856896 |
| 40 | .355496 | .0497871 | .893608 |
| 45 | .580208 | .0497871 | .920880 |
| 50 | .705018 | .0497871 | .941125 |
| 55 | .781435 | .0497871 | .956144 |
| 60 | .831591 | .0497871 | .967298 |
| 70 | .891248 | .0497871 | .981726 |
| 75 | .909828 | .0497871 | .986298 |
| 80 | .924024 | .0497871 | .989706 |
| 85 | .935113 | .0497871 | .992233 |

|        | Exact       | Computed after $B_c(t)$ with $B_c(3) = .0497871$ | Error |
|--------|-------------|--------------------------------------------------|-------|
| $E[B]$ | 19.08553692 | 18.60845683                                      | 3%    |
| $VAR[B]$ | 281.9155718 | 250.9048589                                    | 26%   |

The obtained values, $B_c(t)$, were matched with the lower bounds $B_1(t) = 1 - \frac{e^{2\rho} - 2\rho e^{\rho} - 1}{(1 + \lambda t - e^{\rho})^2}$, remember (4.6), and $B_2(t) = \begin{cases} 0, t < \alpha \\ e^{-\rho}, t \geq \alpha \end{cases}$, see (4.4).

The $B_c(t)$ values always satisfy the ones of $B_1(t)$, sometimes trivial, and of $B_2(t)$.

Note that the busy period of the M|D|∞ queue system has a probability concentration, with value $e^{-\rho}$ in $t = \alpha$, see (5). So, to test the obtained results validity, the busy period mean and variance were computed after the $B_c(t)$ values obtained, but considering $B_c(\alpha) = e^{-\rho}$, and conferred with the real ones.

The values obtained to the mean are very close to the real ones. But the ones obtained to the variance present greater errors. This is natural since the variance calculus accumulates the errors in the computations of the first and second moments centered in the origin. In conclusion, facing the observed errors, it is possible to conclude that results obtained through $B_c(t)$ are quite satisfactory.

Note still that, in principle, the values obtained may be improved decreasing $\Delta t$ (accuracy) and $\Delta p$ (precision). Only in principle, because to run the program is a very slow task and this slowliness increases with the decreasing of $\Delta t$ and $\Delta p$.

### 6. CONCLUSIONS

The M|D|∞ queue busy period is very much studied in the realm of the M|G|∞ queue busy period, of course more general, to work as a reference to it.

As it was seen, it is possible to compute exactly, with no trouble, any parameter related with its distribution, and to perfectly interpret its probability density function.

The only problem stays in the distribution function computation. Precisely, it is presented an algorithm for its computation, which works satisfactorily, apart time problems.

It was also seen that for the traffic intensity great values, above *10*, the M|D|∞ queue busy period is approximately exponentially distributed what simplifies a lot the applicability in practical situations.

### REFERENCES

1. L. Takács. An introduction to queueing theory. Oxford University Press. 1962.
2. M. A. M. Ferreira, M. Andrade. The ties between the M|G|∞ queue system transient behavior and the busy period. International Journal of Academic

Research, 1(1), 84-92, 2009.
3. M. A. M. Ferreira, M. Andrade. M|G|∞ queue system parameters for a

particular collection of service time distributions**.** AJMCSR-African Journal of Mathematics and Computer Science Research, 2(7) 138-141, 2009.

4. M. A. M. Ferreira, M. Andrade. Busy period and busy cycle distributions and parameters for a particular M|G|∞ queue system. American Journal of Mathematics and Statistics, 2(2), 10-15, 2012. DOI: 10.5923/j.ajms.20120202.03.

5. M. A. M. Ferreira. Comportamento transeunte e período de ocupação sistemas de fila de espera sem espera. PhD thesis. ISCTE, 1995.
6. D. J. Daley. Queueing output processes. Advances in Applied Probability, 8, 395-415, 1976.
7. M. F. Ramalhoto, M. A. M. Ferreira. Some further properties of the busy period of an M|G|∞ queue. Central European Journal for Operations Research and Economics, 4(4), 251-278.
8. M. Kendall, A. Stuart. The advanced theory of statistics-Distribution Theory. Charles Griffin and Co., Ltd., London, 1979.
9. B. J. F. Murteira. Probabilidades e estatística. Vol. 1. Editora McGraw-Hill de Portugal, Lisboa,1979.
10. M. A. M. Ferreira, M. Andrade. The M|G|∞ queue busy period exponentiality. APLIMAT- Journal of Applied Mathematics, 4(3), 249-260, 2011.
11. L. K. Platzman, J. C. Ammons, J. J. Bartholdi III. A simple and efficient algorithm to compute tail probabilities from transforms. Operations Research, 36, 137-144, 1988.
12. M. A. M. Ferreira. Distribuição do período de ocupação da fila de espera M|D|∞. Investigação Operacional, 1(16), 43-55, 1996.